\documentclass[11pt,fleqn]{article}
\usepackage{stmaryrd}
\usepackage{amssymb}
\usepackage{bbm}
\usepackage{mathrsfs}
\usepackage{amsfonts}
\usepackage{amsmath,epsf,eufrak}
\begin{document}
\renewcommand{\theequation}{\arabic {section}.\arabic {equation}}
\newcommand{\bi}{\begin{equation}}
\newcommand{\ei}{\end{equation}}
\date{}
\title{  The critical points of a polynomial
\author{Zaizhao Meng}}
\maketitle \baselineskip 24pt
 \begin{center}
 In this paper, we obtain new results on the critical points of a polynomial,
 these results are useful to the Sendov conjecture.
 \end{center}
\section{Introduction}
\setcounter{equation}{0}

Let $a\in (0,1)$,\\
$p(z)=(z-a)\prod\limits_{k=1}^{n-1}(z-z_{k}),\ \ |z_{k}|\leq 1(k=1,...,n-1)$\\
with\\
$p^{\prime}(z)=n\prod\limits_{j=1}^{n-1}(z-\zeta_{j})$.\\
Let $r_{k}=|a-z_{k}|,\ \rho_{j}=|a-\zeta_{j}|$ for
$k,j=1,2,...,n-1$. By relabeling we suppose that \bi
\rho_{1}\leq\rho_{2}\leq ...\leq\rho_{n-1},\  r_{1}\leq r_{2}\leq
...\leq r_{n-1} . \ei
 We have(see [4],[5],[6])
 \bi
 2\rho_{1}\sin(\frac{\pi}{n})\leq r_{k}\leq 1+a,\ \
k=1,2,...,n-1.
 \ei
{\bf Lemma A.} If $0<a<1$ and $\rho_{1}\geq 1$ , then\\
$|p(z)|>1-(1-\lambda)^{n}$ ,\ for
$|z-a|=\lambda\leq\sin(\pi/n)$.\\
This is Lemma 2.16 of [2](see [1]).\\
 {\bf Theorem.} If
$1-(1-|p(0)|)^{\frac{1}{n}}\leq \lambda\leq\sin(\frac{\pi}{n})$ and
$\lambda<a,\ \rho_{1}\geq 1$, then there exists a critical point
$\zeta_{0}=a+\rho_{0}e^{i\theta_{0}}$ such that $ \Re \zeta_{0}\geq
\frac{1}{2}(a-\frac{\lambda(\lambda+2)}{a})$.\\
This theorem improves the previous known results(see [2],[3],[4]).
\section{ Proof of the Theorem }
\setcounter{equation}{0}
 \ We apply Lemma A to conclude that
\\ $|p(z)|>1-(1-\lambda)^{n},\ |z-a|=\lambda$.\\
Since $p(z)$ is univalent in $|z-a|\leq\lambda$, it follows that
there exists a unique point $z_{0}$ with $|z-a|\leq\lambda$ such
that $p(0)=p(z_{0})$. We assume that $\Im z_{0}\geq 0$(if not,
consider $\overline{p(\overline{z})}$). By a variant of the
Grace-Heawood theorem, there exists a critical point in each of the
half-planes bounded by the perpendicular bisector $L$ of the segment
from 0 to $z_{0}$. Let $ \zeta_{0}=a+\rho_{0}e^{i\theta_{0}}$ be the
critical point in the half-plane containing $z_{0}$.\\
The equation of $L$ is $|z|=|z-z_{0}|$, that is
 \bi
z\overline{z}=(z-z_{0})(\overline{z}-\overline{z_{0}}),
 \ei
 then
 \bi
 |z_{0}|^{2}=z \overline{z_{0}}+\overline{z}z_{0}.
\ei Let $z^{\ast}=e^{i\beta_{0}}$ be the joint point of $L$ and the
circle $|z|=1,\ \Im z^{\ast}\geq 0$, then

$|z_{0}|^{2}=e^{i\beta_{0}}\overline{z_{0}}+e^{-i\beta_{0}}z_{0}$.

Hence

 $e^{i\beta_{0}}=\frac{z_{0}}{2}\pm\frac{z_{0}}{2|z_{0}|}\sqrt{|z_{0}|^{2}-4}$,\\
 that is
 \bi
e^{i\beta_{0}}=(\frac{1}{2}\pm\frac{i}{2|z_{0}|}\sqrt{4-|z_{0}|^{2}})z_{0}.
 \ei

Write $z_{0}=a+re^{i\alpha},\ r\leq\lambda$, then we choose\\

$
\cos\beta_{0}=\frac{1}{2}(a+r\cos\alpha)-\frac{1}{2|z_{0}|}\sqrt{4-|z_{0}|^{2}}
r\sin\alpha, $

\bi
\cos\beta_{0}=\frac{1}{2}(a+r\cos\alpha)-\frac{1}{2}\frac{\sqrt{4-a^{2}-2ar\cos\alpha-r^{2}}}{\sqrt{a^{2}+2ar\cos\alpha+r^{2}}}r\sin\alpha,
\ei

We fix $r$ and consider the circle $|z-a|=r$, let
$\sin\alpha_{1}=\frac{\sqrt{a^{2}-r^{2}}}{a} , \ x=\cos\alpha$ and $
F(x)=x-\sqrt{\frac{4-a^{2}-2arx-r^{2}}{a^{2}+2arx+r^{2}}}\sqrt{1-x^{2}}$,\\
then
 \bi
 \cos\beta_{0}=\frac{1}{2}(a+rF(x)),
 \ei
it is sufficient to give lower bound of $F(x)$ for
$x\in[-1,-\frac{r}{a}]$.

Write $G(x)=-F(-x)$, then
 \bi
G(x)=x+\sqrt{\frac{4-a^{2}+2arx-r^{2}}{a^{2}-2arx+r^{2}}}\sqrt{1-x^{2}},
\ei

it is sufficient to give upper bound of $G(x)$ for $x\in
[\frac{r}{a},1]$.\\
Write $\psi=\frac{4-a^{2}+2arx-r^{2}}{a^{2}-2arx+r^{2}},\
\phi=a^{2}-2arx+r^{2}.$\\
We have
 \bi
G^{\prime}(x)=1+\frac{1}{2}\psi^{-\frac{1}{2}}\psi^{\prime}\sqrt{1-x^{2}}-\psi^{\frac{1}{2}}x(1-x^{2})^{-\frac{1}{2}},
\ei

and $\psi^{\prime}=8ar(a^{2}-2arx+r^{2})^{-2}$,\\
hence
\bi
 G^{\prime}(x)=\phi^{-2}\psi^{-\frac{1}{2}}(1-x^{2})^{-\frac{1}{2}}G_{1},
 \ei

where

\bi
G_{1}=\phi^{\frac{3}{2}}(4-\phi)^{\frac{1}{2}}(1-x^{2})^{\frac{1}{2}}+4ar(1-x^{2})-x\phi(4-\phi).
\ei

We have $x=\frac{a^{2}+r^{2}-\phi}{2ar}$, then

\bi
 G_{1}=\frac{1}{2ar}G_{2},
\ei

where\\

$
G_{2}=\phi^{\frac{3}{2}}(4-\phi)^{\frac{1}{2}}(4a^{2}r^{2}-(a^{2}+r^{2})^{2}+2(a^{2}+r^{2})\phi-\phi^{2})^{\frac{1}{2}}\\
+8a^{2}r^{2}-2(a^{2}+r^{2})^{2}+4(a^{2}+r^{2})\phi-2\phi^{2}-(a^{2}+r^{2}-\phi)\phi(4-\phi)
$, \\ hence\\
$G_{2}=\phi^{\frac{3}{2}}(4-\phi)^{\frac{1}{2}}(4a^{2}r^{2}-(a^{2}+r^{2})^{2}+2(a^{2}+r^{2})\phi-\phi^{2})^{\frac{1}{2}}\\
+8a^{2}r^{2}-2(a^{2}+r^{2})^{2}+(a^{2}+r^{2}+2)\phi^{2}-\phi^{3},
$\\
and $\phi\in [(a-r)^{2},a^{2}-r^{2}].$\\
The roots of $G_{2}$ satisfy\\
$
L=\phi^{3}(4-\phi)(4a^{2}r^{2}-(a^{2}+r^{2})^{2}+2(a^{2}+r^{2})\phi-\phi^{2})=\\
(\phi^{3}-(a^{2}+r^{2}+2)\phi^{2}+2(a^{2}+r^{2})^{2}-8a^{2}r^{2})^{2}=R,
$ \\
say.\\
We have\\
$L=\phi^{6}-2(a^{2}+r^{2}+2)\phi^{5}+(8(a^{2}+r^{2})-4a^{2}r^{2}+(a^{2}+r^{2})^{2})\phi^{4}+
(16a^{2}r^{2}-4(a^{2}+r^{2})^{2})\phi^{3}$,\\
$R=\phi^{6}-2(a^{2}+r^{2}+2)\phi^{5}+(a^{2}+r^{2}+2)^{2}\phi^{4}+4(a^{2}-r^{2})^{2}\phi^{3}
-4(a^{2}+r^{2}+2)(a^{2}-r^{2})^{2}\phi^{2}+4(a^{2}-r^{2})^{4}$.\\
Let $d=1-r^{2}, \ e=a^{2}-1$, by $L=R$, we deduce\\
$ed\phi^{4}-2(e+d)^{2}\phi^{3}+(4+e-d)(e+d)^{2}\phi^{2}-(e+d)^{4}=0$,\\
that is\\
$(e\phi^{2}-2(e+d)\phi-(e+d)^{2})(d\phi^{2}-2(e+d)\phi+(e+d)^{2})=0$,\\
there is only the root $\phi_{0}=\frac{a^{2}-r^{2}}{1+r}\in
[(a-r)^{2},a^{2}-r^{2}]$, and\\
$x_{0}=\frac{a^{2}+r^{2}-\phi_{0}}{2ar}=\frac{2r+a^{2}+r^{2}}{2a(1+r)}\in
[\frac{r}{a},1]$.\\
We have $G^{\prime}(1)<0,\ G^{\prime}(\frac{r}{a})>0$, $G(x_{0})$ is
the maxima value of $G(x)$ for $x\in [\frac{r}{a},1]$.\\
We obtain  $G(x_{0})=\frac{r+2}{a}$,and $G(x)\leq\frac{r+2}{a}, \
x\in [\frac{r}{a},1]$,
 hence\\
 $F(x)\geq -\frac{r+2}{a},\ x\in [-1,-\frac{r}{a}],\\
\cos\beta_{0}\geq\frac{1}{2}(a-\frac{r(r+2)}{a})\geq\frac{1}{2}(a-\frac{\lambda(\lambda+2)}{a}),$\\
the theorem follows.\\

{ {\small E-mail:\ mengzzh@126.com}
\end{document}